\newcommand{\tikzcuboid}[5]{
\begin{tikzpicture}[scale=#4]
\foreach \x in {0,...,#1}
{   \draw (\x ,0  ,#3 ) -- (\x ,#2 ,#3 );
    \draw (\x ,#2 ,#3 ) -- (\x ,#2 ,0  );
}
\foreach \x in {0,...,#2}
{   \draw (#1 ,\x ,#3 ) -- (#1 ,\x ,0  );
    \draw (0  ,\x ,#3 ) -- (#1 ,\x ,#3 );
}
\foreach \x in {0,...,#3}
{   \draw (#1 ,0  ,\x ) -- (#1 ,#2 ,\x );
    \draw (0  ,#2 ,\x ) -- (#1 ,#2 ,\x );
}

\foreach \x in {0,...,#1}{
    \pgfmathtruncatemacro{\j}{mod(\x,#5)}
    \foreach \y in {0,...,#2}{
        \ifthenelse{\equal{\j}{0}}
        {\filldraw[red] (\x,\y,#3) circle (4pt);}
        {\filldraw[blue] (\x,\y,#3) circle (4pt);}
    }
}
\foreach \y in {0,...,#2}{
    \foreach \z in {0,...,#3}{
        \filldraw[red] (#1,\y,\z) circle (4pt);
    }
}
\foreach \x in {0,...,#1}{
    \pgfmathtruncatemacro{\j}{mod(\x,#5)}
    \foreach \z in {0,...,#3}{
        \ifthenelse{\equal{\j}{0}}
        {\filldraw[red] (\x,#2,\z) circle (4pt);}
        {\filldraw[blue] (\x,#2,\z) circle (4pt);}
    }
}
\pgfmathtruncatemacro{\numbuf}{#1/#5}
\pgfmathsetmacro{\halfbuf}{#5/2}
\foreach \x in {0,...,\numbuf}{
    \pgfmathtruncatemacro{\odd}{2*\x+1}
    \node [below,red] at (\x*#5,0,#3) {$I_{\odd}$};
}
\pgfmathtruncatemacro{\numbuf}{#1/#5 - 1}
\pgfmathsetmacro{\halfbuf}{#5/2}
\foreach \x in {0,...,\numbuf}{
    \pgfmathtruncatemacro{\even}{2*\x+2}
    \node [below,blue] at (\x*#5+\halfbuf,0,#3) {$I_{\even}$};
}

\draw [pen colour={black},
    decorate,
    thick,
    decoration = {calligraphic brace,
        raise=5pt,
        amplitude=5pt,
        aspect=0.5,mirror}] (#1,0) --  (#1,#2)
node[pos=0.5,right=10pt,black]{$n$};

\draw [pen colour={black},
    decorate,
    thick,
    decoration = {calligraphic brace,
        raise=5pt,
        amplitude=5pt,
        aspect=0.5}] (#1-#5,#2) --  (#1,#2)
node[pos=0.5,above=10pt,black]{$b$};
\end{tikzpicture}
}

\documentclass[graybox]{svmult}

\usepackage{type1cm}        
\usepackage{makeidx}         
\usepackage{graphicx}        
\usepackage{multicol}        
\usepackage[bottom]{footmisc}
\usepackage{hyperref}

\usepackage{newtxtext}       %
\usepackage[varvw]{newtxmath}       
\usepackage{amsmath,bm,sidecap}
\usepackage{tikz}
\usetikzlibrary{decorations.pathreplacing,calligraphy}
\usetikzlibrary{calc,patterns}
\usepackage[linesnumbered,ruled,vlined]{algorithm2e}

\makeindex             

\newcommand{\mtx}[1]{\bm{\mathsf{#1}}}
\newcommand{\vct}[1]{\bm{\mathsf{#1}}}

\begin{document}

\title*{GPU Optimizations for the Hierarchical Poincaré-Steklov Scheme}
\author{Anna Yesypenko and Per-Gunnar Martinsson}
\institute{Anna Yesypenko \at Oden Institute, \email{annayesy@utexas.edu}
\and Per-Gunnar Martinsson \at Oden Institute \email{pgm@oden.utexas.edu}}
%
%
\maketitle

\abstract{This manuscript presents GPU optimizations for the 2D Hierarchical Poincaré-Steklov (HPS) discretization scheme. HPS is a multi-domain spectral collocation method that combines high-order discretizations with direct solvers to accurately resolve highly oscillatory solutions. The domain decomposition approach of HPS connects domains directly via a sparse direct solver. The proposed optimizations exploit batched linear algebra on modern hybrid architectures, are straightforward to implement, and improve the solver's practical speed. The manuscript demonstrates that GPU optimizations can significantly reduce the traditionally high cost of performing local static condensation for discretizations with very high local order $p$. Numerical experiments for the Helmholtz equation with high wavenumbers on curved and rectangular domains confirm the high accuracy achieved by the HPS discretization and the significant reduction in computation time achieved with GPU optimizations.
}

\abstract*{This manuscript presents GPU optimizations for the 2D Hierarchical Poincaré-Steklov (HPS) discretization scheme. HPS is a multi-domain spectral collocation method that combines high-order discretizations with direct solvers to accurately resolve highly oscillatory solutions. The domain decomposition approach of HPS connects domains directly via a sparse direct solver. The proposed optimizations exploit batched linear algebra on modern hybrid architectures, are straightforward to implement, and improve the solver's practical speed. The manuscript demonstrates that GPU optimizations can significantly reduce the traditionally high cost of performing local static condensation for discretizations with very high local order $p$. Numerical experiments for the Helmholtz equation with high wavenumbers on curved and rectangular domains confirm the high accuracy achieved by the HPS discretization and the significant reduction in computation time achieved with GPU optimizations.
}

\section{Introduction}
\label{sec:intro}

We describe methods for solving boundary value problems of the form
\begin{equation}
\left\{\begin{aligned}
\mathcal A u(x) =&\ f(x),\qquad&x \in \Omega,\\
u(x) =&\ g(x),\qquad&x \in \partial \Omega,
\end{aligned}\right.
\label{eq:bvp}
\end{equation}
where $\mathcal A$ is a second order elliptic differential operator,
and $\Omega$ is domain in two dimensions with boundary $\partial \Omega$.
For the sake of concreteness, we will focus on the case where $\mathcal A$ is a variable
coefficient Helmholtz operator
\begin{equation}
\label{eq:var_helm}
\mathcal A u(x) = -\Delta u(x) - \kappa^{2}b(x)u(x),
\end{equation}
where $\kappa$ is a reference wavenumber,
and where $b(x)$ is a smooth non-negative function that typically satisfies $0 \le b(x) \le 1$.
Upon discretizing (\ref{eq:bvp}), one obtains a linear system
\begin{equation}
\label{eq:Au=f}
\mtx { A}\vct{u} = \vct{f}
\end{equation}
involving a sparse coefficient matrix
$\mtx {A} \in \mathbb{R}^{N \times N}$.
The focus of this work is on efficiently solving the sparse system
(\ref{eq:Au=f}) for the Hierarchical Poincar\'e-Steklov (HPS) discretization.
HPS is a multi-domain spectral collocation scheme that allows for relatively
high choices of $p$, while interfacing well with sparse direct solvers.
For (\ref{eq:bvp}) discretized with HPS with local polynomial order $p$, 
the cost of factorizing $\mtx { A}$ directly is 
\begin{equation}
T_{\rm build} = O \left ( \underset{\rm leaf\ operations} {p^4 N} + \underset{\rm direct\ solver}{N^{3/2}} \right).
\end{equation}
After the leaf operations are complete, the cost to factorize the system directly has no pre-factor dependence on $p$.
The pre-factor cost of the leaf operations, however, has long been viewed as prohibitively expensive.
This manuscript describes simple GPU optimizations using batched linear algebra that substantially 
accelerate the leaf operations and shows compelling results for $p$ up to 42.
We also demonstrate that the choice of $p$ does not have substantial effects on the build time for the direct factorization stage,
allowing $p$ to be chosen based on physical considerations instead of practical concerns.

High order discretization is crucial in resolving variable-coefficient
scattering phenomena due to the well known ``pollution effect'' that generally requires the number of points per wavelength to increase, the larger the computational domain is.
The pollution effect is very strong for low order discretizations, but quickly gets less problematic as the discretization order increases \cite{beriot2016efficient,deraemaeker1999dispersion}.
HPS is less sensitive to pollution because the scheme allows for high choices of local polynomial order $p$ \cite{gillman2014direct,martinsson2013direct}.
Combining HPS discretization with efficient sparse direct solvers provides a powerful tool for 
resolving challenging scattering phenomena to high accuracy, especially for situations where
no efficient preconditioners are known to exist (e.g. trapped rays, multiple reflections, backscattering) \cite{ernst2012difficult}.

\section{HPS Discretization and Interfacing with Sparse Direct Solvers}
\label{sec:disc}

We next discuss the HPS discretization and efficient methods to interface the resulting
sparse linear system with direct solvers.
We introduce the HPS briefly for the simple model problem (\ref{eq:bvp}), 
and refer the reader to \cite{babb2018accelerated,2019_martinsson_book,hao2016direct} for details and extensions. 
An important limitation of the discretization is that we assume the solution is smooth
and that the coefficients in the operator $\mathcal A$ of (\ref{eq:bvp}) are smooth as well.

\begin{figure}[htb!]
\centering
\begin{minipage}{3.5cm}
\centering
\includegraphics[width=3cm]{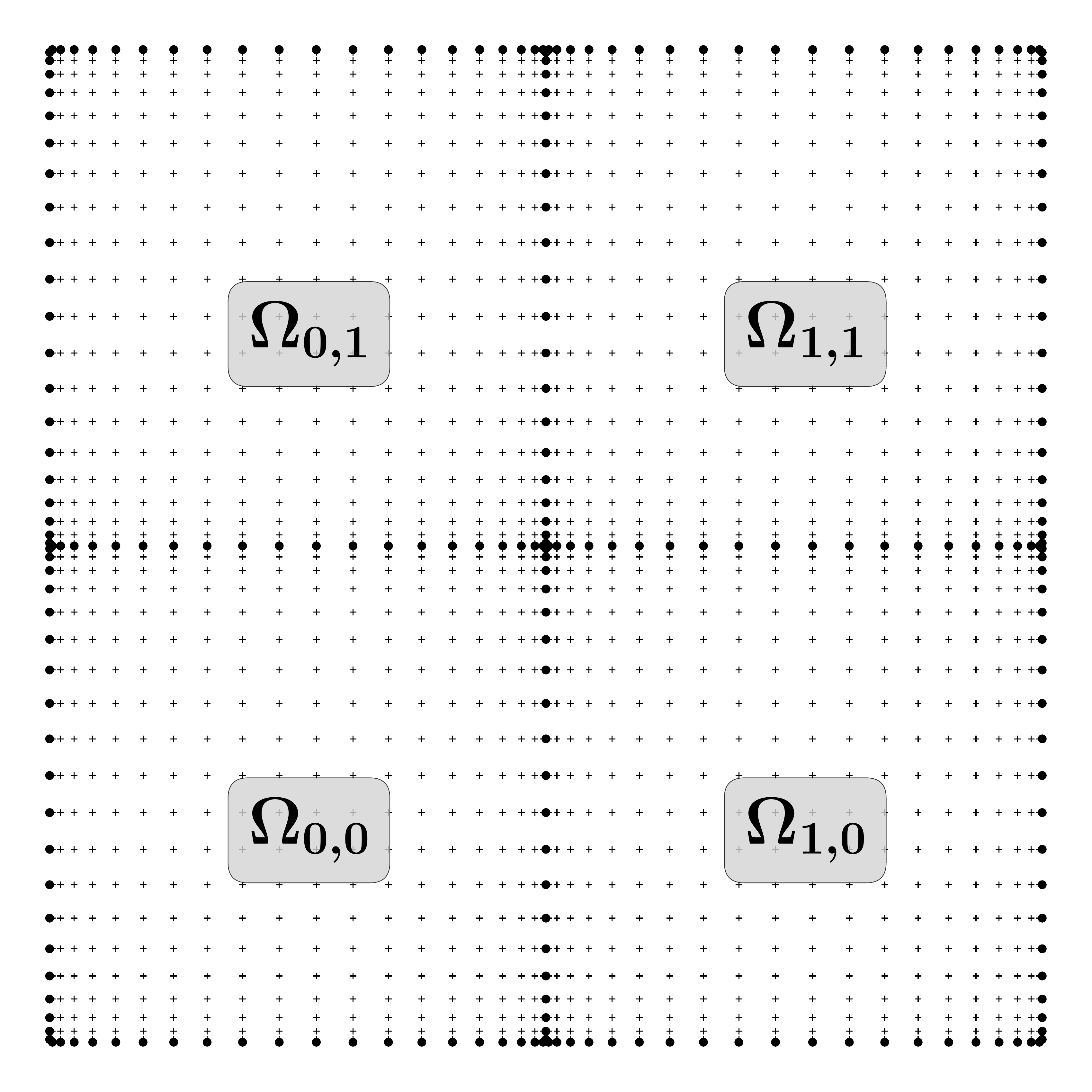}
\end{minipage}%
$\Rightarrow$
\begin{minipage}{3.5cm}
\centering
\includegraphics[width=3cm]{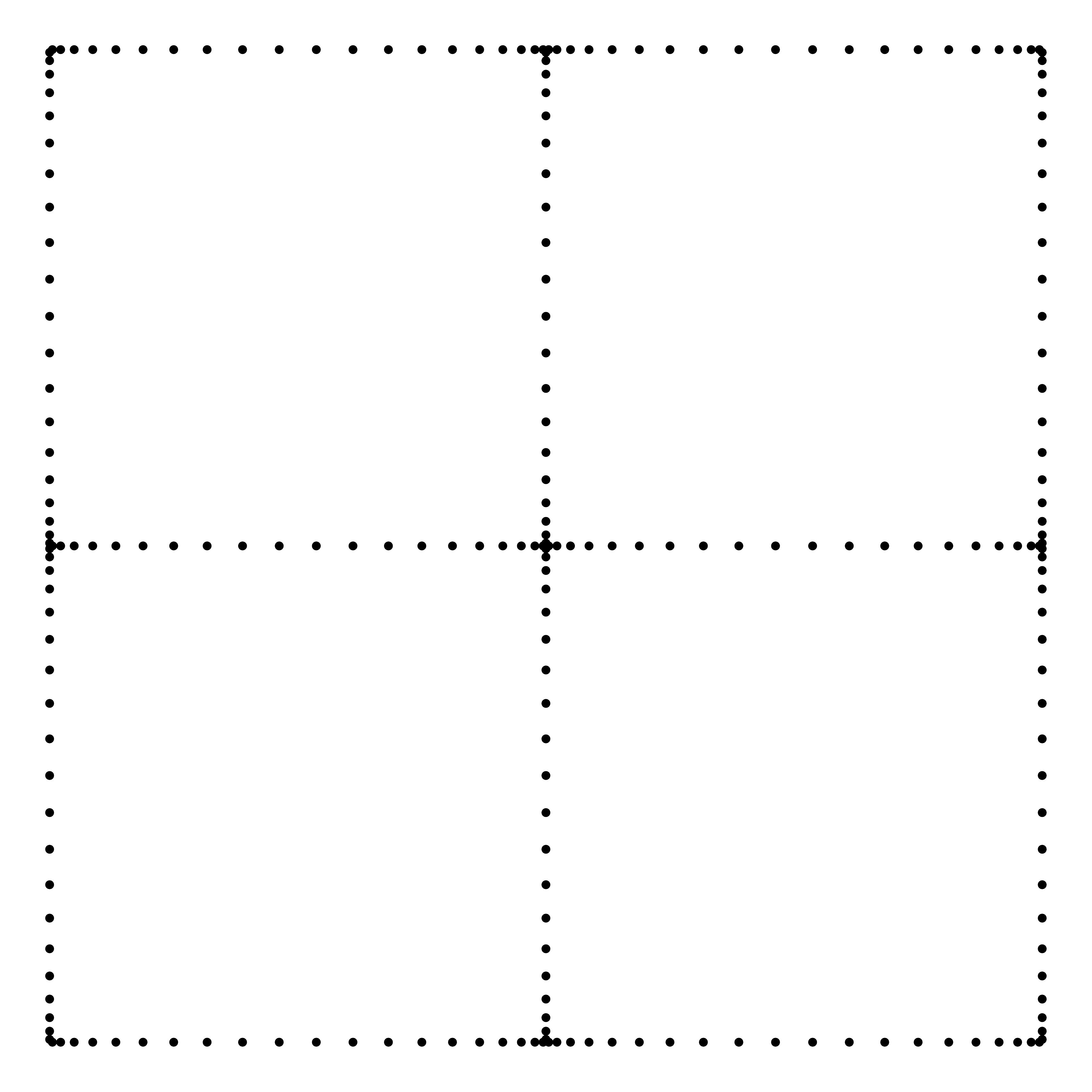}
\end{minipage}%
\hfill
\begin{minipage}{4cm}
\caption{\small Prior to interfacing with sparse direct solvers, we do static condensation to produce an equivalent system (\ref{eq:hps_reduced}) to solve on the remaining active nodes. The original grid has $N$ points, and remaining grid has $\approx N/p$ points.}
\end{minipage}
\label{fig:hps_reduced_grid}
\end{figure}

The domain $\Omega$ is partitioned into non-overlapping subdomains. 
The discretization is described by two parameters, $a$ and $p$, which are the element size and local
polynomial order, respectively.
On each subdomain, we place a $p \times p$
tensor product mesh of Chebyshev points. Internal to each subdomain, the PDE is enforced locally via
spectral differentiation and direct collocation. On element boundaries, we enforce that the flux between adjacent
boundaries is continuous.
On each subdomain of $p^2$ nodes, the spectral differentiation operators lead to a dense
matrix of interactions of size $p^2 \times p^2$.
To improve efficiency of sparse direct solvers for HPS discretizations, we
``eliminate'' the dense interactions of nodes interior to each subdomain. This process is referred to as ``static condensation'' \cite{cockburn2016static,guyan1965reduction}.
The remaining active nodes are on the boundaries between subdomains. 
As a result of the leaf elimination,
we produce a smaller system $\tilde {\mtx A}$ of size $\approx N / p$ with equivalent body load
$\tilde {\mtx f}$ on the active nodes located on the boundaries between subdomains, as shown in Figure \ref{fig:hps_reduced_grid}
\begin{equation}\tilde {\mtx A} \tilde {\mtx u} = \tilde {\mtx f}.\label{eq:hps_reduced}\end{equation}

Due to the domain decomposition used in HPS, the leaf operations required to produce the equivalent system (\ref{eq:hps_reduced}) can be done embarrassingly in parallel. The leaf operations require independent dense linear algebraic operations (e.g., LU factorization, matrix-matrix multiply) on $N / {p^2}$ systems, each of size $p^2 \times p^2$, resulting in an overall cost of $O(p^4 N)$. For $p$ up to about 42, these operations can be efficiently parallelized with batched linear algebra (BLAS). However, for larger $p$, methods that produce a sparser equivalent system may be more appropriate \cite{brubeck2021scalable,2021_fortunato_ultraspherical}.

Overhead costs can make achieving high arithmetic intensity for many small parallel tasks a challenge. However, batched BLAS offers a solution. It is highly optimized software for parallel operations on matrices that are small enough to fit in the top levels of the memory hierarchy (i.e., smaller than the L2 cache size) \cite{dongarra2017design}. The framework groups small inputs into larger "batches" to automatically achieve good parallel performance on high-throughput architectures such as GPUs.

The technique we present is most readily applicable to the case where the same discretization order $p$ is used on every discretization patch. However, it would not be too difficult to allow $p$ to be chosen from a fixed set of values (say $ p \in \{6,10,18,36\}$ or something similar). This would enable many of the advantages of hp-adaptivity, while still enabling batching to accelerate computations. 

\begin{remark}
Since the leaf computations are \textit{very} efficient, we saved memory and reduced communication by not explicitly storing the factorizations of the local spectral differentiation matrices.
Instead, these are reformed and refactored after each solve involving the reduced system  (\ref{eq:hps_reduced}).
\end{remark}

We combined the fast leaf factorization procedure with two methods for solving the reduced system (\ref{eq:hps_reduced}). 
The first option for solving (\ref{eq:hps_reduced}) uses a black-box sparse direct solver with 
the nested dissection (ND) ordering.
ND is a based on a multi-level graph partitioning of nodes and produces a sparse factorization with 
minimal fill-in \cite{2006_davis_directsolverbook,2016_acta_sparse_direct_survey}.
In 2D, sparse factorization using the ND ordering requires $O \left( N^{3/2} \right)$ time to build and $O \left( N \log N \right)$ time to apply.

\begin{figure}[htb!]
\centering
\begin{minipage}{6.5cm}
\centering
\tikzcuboid{16}{7}{0}{0.3}{4}
\end{minipage}%
\hfill
\begin{minipage}{5.0cm}
\caption{\small Domain decomposition used in SlabLU. The even-numbered nodes correspond to the nodes interior to each subdomain. The odd-numbered nodes correspond to interfaces between slabs. The slab partitioning is chosen so that interactions between slab interiors are zero. The slabs have width of $b$ points.}
\end{minipage}
\label{fig:slablu_decomp}
\end{figure}

As a second option for solving (\ref{eq:hps_reduced}), we used a 
scheme we refer to as SlabLU, which is a simplified two-level scheme
(as opposed to standard hierarchical schemes) that is designed for
ease of parallelization \cite{yesypenko2022slablu}.
To be precise, SlabLU uses a decomposition of the domain into elongated ``slab'' subdomains, as shown in Figure \ref{fig:slablu_decomp}. With this decomposition, the linear system (\ref{eq:hps_reduced}) has the block form 
\begin{equation}
\label{eq:dell1}
\begin{bmatrix}
\tilde {\mtx A}_{11} & \tilde {\mtx A}_{12} & \mtx 0 & \mtx 0 & \mtx 0 & \dots \\
\tilde {\mtx A}_{21} & \tilde {\mtx A}_{22} & \tilde {\mtx A}_{23} & \mtx 0 &  \mtx 0 & \dots \\
\mtx 0 & \tilde {\mtx A}_{32} & \tilde {\mtx A}_{33} & \tilde {\mtx A}_{34} &   \mtx 0 & \dots\\
\mtx{0} & \mtx 0 & \tilde {\mtx A}_{43} & \tilde {\mtx A}_{44} & \tilde {\mtx A}_{45} & \dots\\
\vdots & \vdots & \vdots & \vdots & \vdots & \vdots
\end{bmatrix}
\begin{bmatrix} \tilde {\mtx u}_1\\\tilde {\mtx u}_2\\ \tilde {\mtx u}_3 \\ \tilde {\mtx u}_4 \\ \vdots \end{bmatrix} =
\begin{bmatrix} \tilde {\mtx f}_1\\ \tilde {\mtx f}_2\\ \tilde {\mtx f}_3\\ \tilde {\mtx f}_4 \\ \vdots        \end{bmatrix}.
\end{equation}
The nodes internal to each slab are eliminated by computing sparse factorizations 
of the diagonal blocks $\tilde {\mtx A}_{22}, \tilde {\mtx A}_{44}, \dots$ in parallel. This results in another block tridiagonal coefficient matrix $\mtx{T}$ that has much smaller blocks than $\tilde{\mtx{A}}$ (and half as many). The blocks of $\mtx{T}$ are dense, but can be represented efficiently using data sparse formats such as the $\mathcal{H}$-matrix format of Hackbusch. The ranks are very low, due to the thinness of the slabs. The construction of these blocks is further accelerated by using the black box randomized compression techniques described in  
\cite{levitt2022linear}.

The reduced linear system with blocks having $\mathcal H$-matrix structure can in principle be solved efficiently using rank-structured linear algebra. However, we found that for 2D problems, it is most efficient to relinquish the rank structure and simply convert all blocks to a dense format before factorizing the block tridiagonal system.
(In 3D, this simplistic approach is possible only for small problems.)
With a choice of slab width $b$ that grows slowly with the problem size as $b \sim n^{2/3}$, the resulting two-level scheme has complexity
$O(N^{5/3})$ to factorize ${\tilde {\mtx A}}$ directly and $O \left( N^{7/6} \right)$
complexity to apply the computed factors to solve (\ref{eq:hps_reduced}). SlabLU is simple scheme that leverages high concurrency and
batched BLAS to achieve high performance on modern hybrid architectures.
Despite the asymptotically higher costs, SlabLU performs favorably compared to multi-level
nested dissection schemes in its build time and memory footprint, as we demonstrate in Section \ref{sec:numerical}.
\cite{yesypenko2022slablu} provides details on SlabLU.

\section{Numerical Experiments}
\label{sec:numerical}

We demonstrate the effectiveness of the HPS discretization combined with
sparse direct solvers in solving high-frequency Helmholtz equations. 
The experiments were conducted on a desktop computer equipped with a 16-core Intel i9-12900k CPU
and 128GB of memory, and a NVIDIA RTX 3090 GPU with 24GB of memory.

\begin{figure}[htb!]
\centering
\begin{minipage}{7.8cm}
\includegraphics[width=\textwidth]{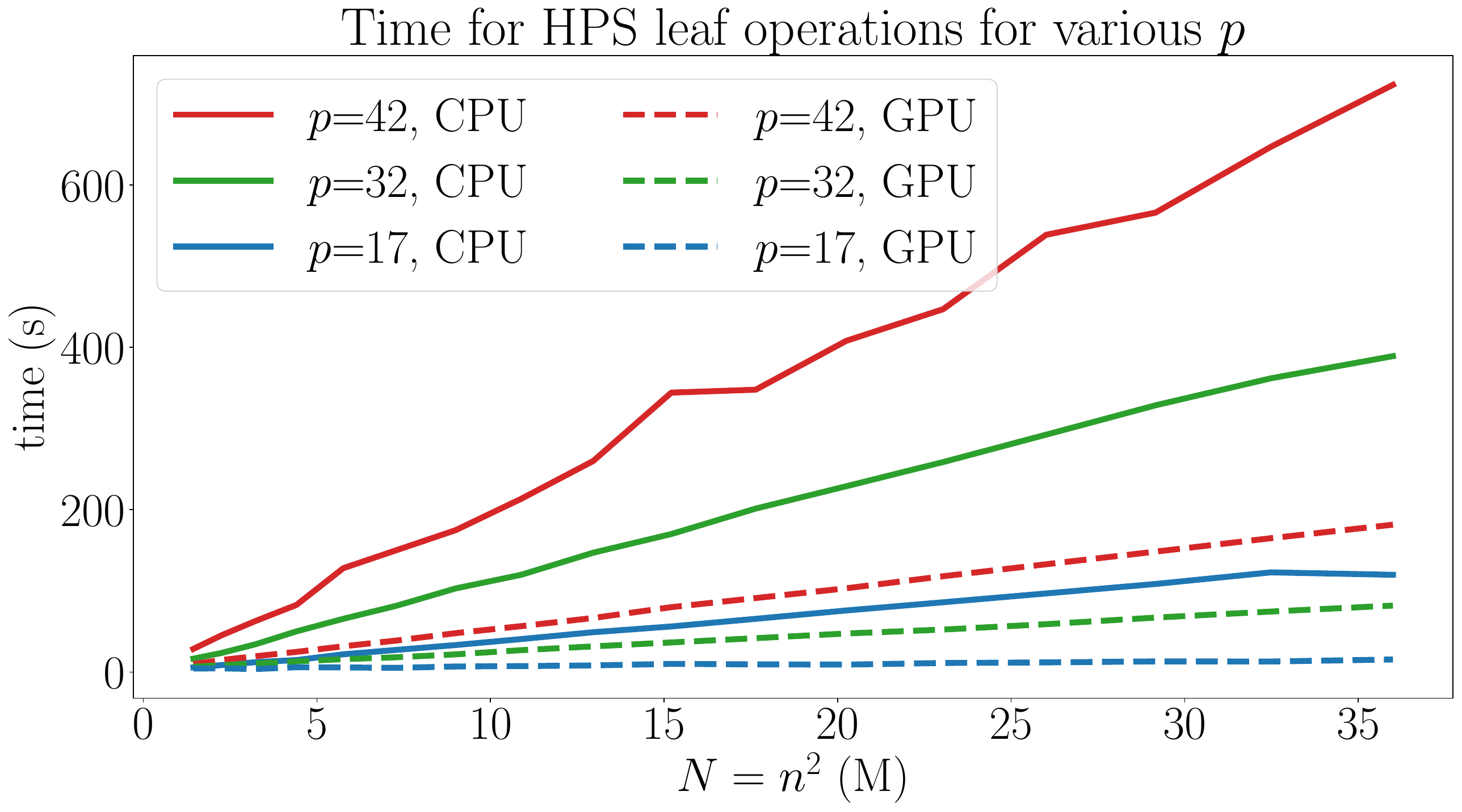}
\end{minipage}%
\hfill
\begin{minipage}{3.5cm}
\caption{Leaf operations for HPS require $O(p^4 N)$ operations, though the practical scaling
for parallel operations
has a small constant prefactor for $p$ up to 42.
Parallel HPS leaf operations are further accelerated on GPUs, with a speed-up of at least 4x.}
\label{fig:hps_batchedblas}
\end{minipage}
\end{figure}

We show that GPU optimizations enable efficient leaf operations for various local polynomial orders,
cf.~Figure \ref{fig:hps_batchedblas}. After the leaf operations, we directly factorize the reduced
system (\ref{eq:hps_reduced}) using efficient sparse direct solvers. We demonstrate that the choice of $p$
does not significantly affect the time to factorize ${\tilde {\mtx A}}$. Having the freedom to choose $p$
allows the user to resolve highly oscillatory PDEs to high-order accuracy without worrying about how 
the choice may affect the cost of solving (\ref{eq:Au=f}) directly.

To demonstrate the effectiveness of the HPS discretization resolving
oscillatory solutions to high accuracy, we report results for a PDE with a known analytic solution
\begin{equation}
\left\{\begin{aligned}
-\Delta u(x)  - \kappa^2 u(x) =&\ 0,\qquad&x \in \Omega={[0,1]}^2,\\
u(x) =&\ u_{\text{true}}(x),\qquad&x \in \partial \Omega,
\end{aligned}\right. 
\label{eq:constant_coeff_helm}
\end{equation}
The true solution $u_{\text{true}}$ is given by 
$u_{\text{true}} = J_0 \left( \kappa |x - (-0.1,0.5)| \right)$,
where $x \mapsto J_{0}(\kappa|x|)$ is the free-space fundamental solution
to the Helmholtz equation.  We discretize (\ref{eq:constant_coeff_helm}) using HPS for various choices
of $p$ and set the wavenumber $\kappa$ to increase with $N$ to maintain 10 points per wavelength with 
increasing problem size.
After applying a direct solver to solve (\ref{eq:hps_reduced}) on the reduced HPS grid, 
we re-factorize the linear systems on interior leaf nodes to calculate the solution $\mtx u_{\rm calc}$ 
on the full HPS grid. The leaf solve requires time $O(p^4 N)$ but is particularly efficient using the 
GPU optimizations described. The reported build times and solve times include the leaf operations. 
We report the relative error with respect to the residual of the discretized system (\ref{eq:Au=f}). 
When a true solution is known, we also report the relative error with respect to the true 
solution $\mtx u_{\rm true}$ evaluated on the collocation points of the full HPS grid
\begin{equation} 
\text{relerr}_{\text{res}} = \frac{{\|\mtx { A} \mtx u_{\text{calc}} - \mtx f\|}_2}{{\|\mtx f\|}_2}, \qquad \text{relerr}_{ \rm true} = \frac{{\|\mtx u_{\rm calc} - \mtx u_{\rm true}\|}_2}{{\| \mtx u_{\rm true}\|}_2}.
\label{eq:accuracy_formula}
\end{equation}

\subsection{Comparison of Sparse Direct Solvers}

The system (\ref{eq:hps_reduced}) is solved using two different sparse direct solvers, SuperLU and SlabLU.
SuperLU is a black-box solver that finds an appropriate ordering of the system to minimize fill-in while
increasing concurrency by grouping nodes into super-nodes \cite{2011_li_supernodal}.
We accessed SuperLU through the Scipy interface (version 1.8.1) and called it with the COLAMD ordering.
This version of Scipy uses the CPU only. Not many GPU-aware sparse direct solvers are widely available,
though this is an active area of research. 
SuperLU uses a pivoting scheme that can exchange rows between super-nodes
to attain almost machine precision accuracy in the residual of the computed solutions.

\begin{figure}[htb!]
\begin{minipage}{7.3cm}
\includegraphics[width=7.3cm]{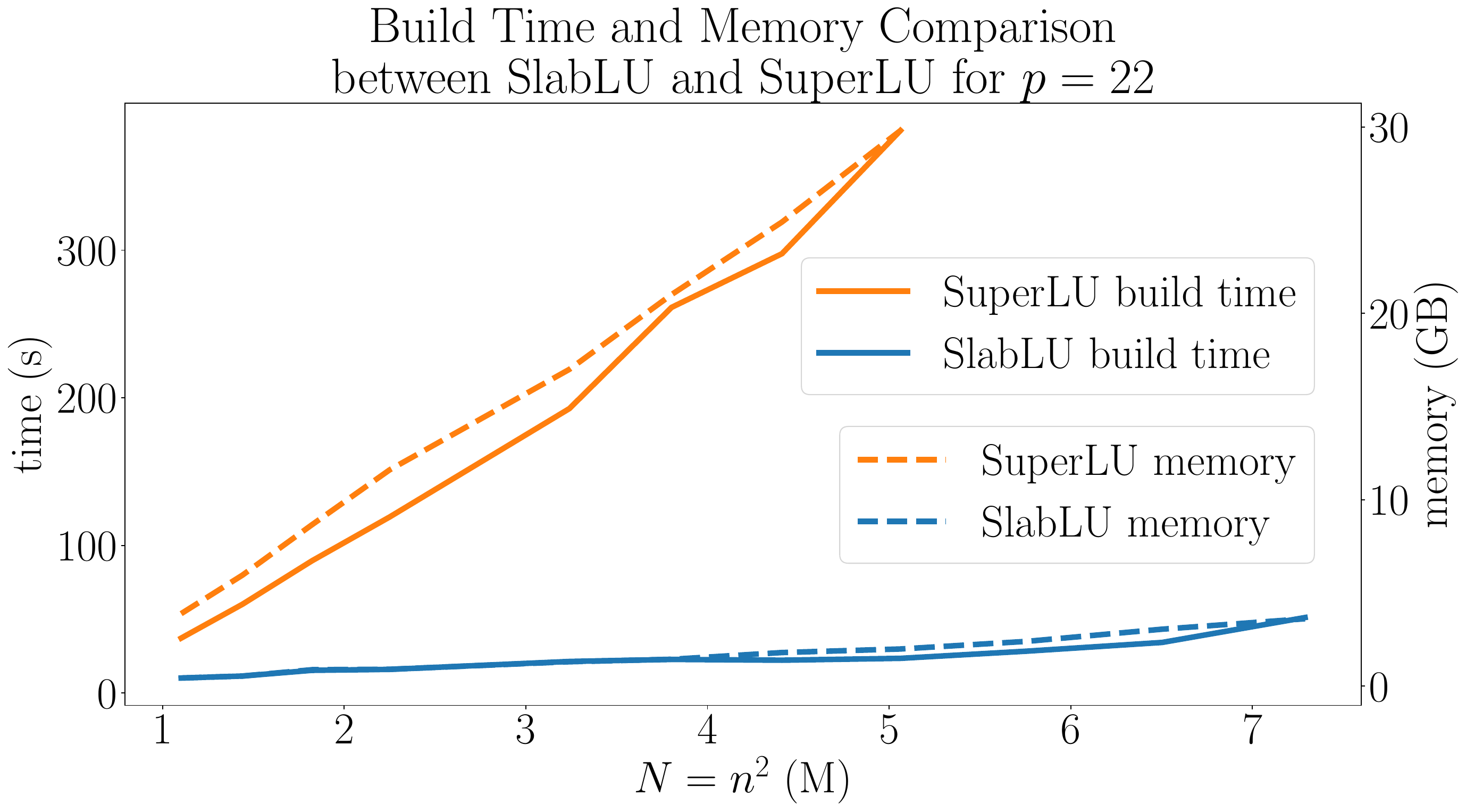}
\end{minipage}%
\hfill
\begin{minipage}{4.1cm}
\caption{Build time and memory footprint comparison between SuperLU and SlabLU for (\ref{eq:constant_coeff_helm}) 
discretized with HPS ($p$=22), where $\kappa$ is increased
to maintain 10 points per wavelength.
For $N$=5.06M, the SlabLU build time is faster by a factor of 16x
and the memory footprint is less by a factor of 14x. }
\label{fig:comparison_build}
\end{minipage}
\vspace{1em}

\includegraphics[width=0.49\textwidth]{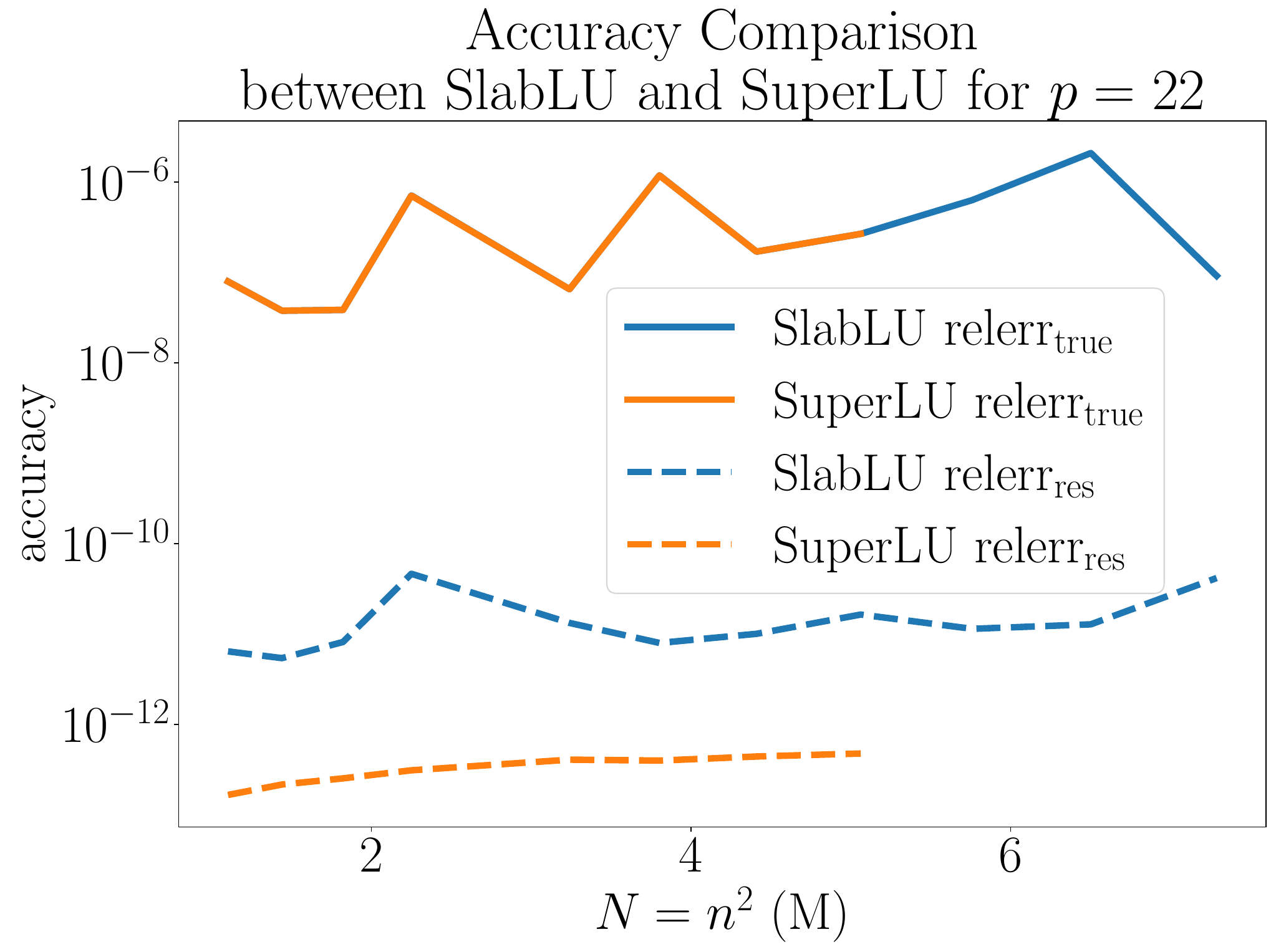}
\includegraphics[width=0.49\textwidth]{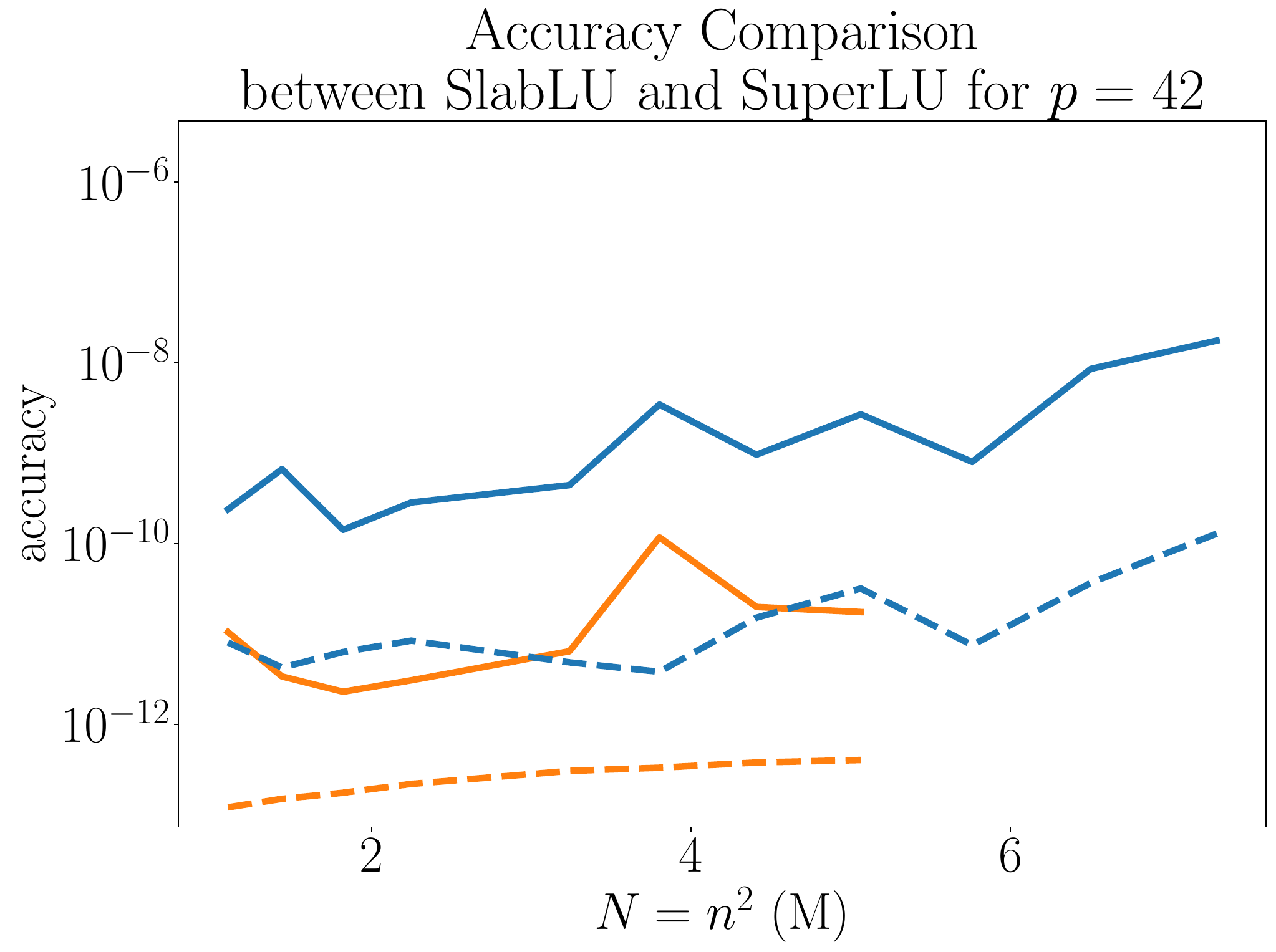}
\caption{Accuracy comparison between SuperLU and SlabLU for (\ref{eq:constant_coeff_helm}) 
discretized with HPS ($p$=22,42), where $\kappa$ is increased
to maintain 10 points per wavelength. SuperLU uses a more sophisticated pivoting scheme to
achieve high accuracy in the residual error. For $p$=22, SlabLU and SuperLU both resolve 
the solution to 6 digits in the true relative accuracy.
For $p$=42, SuperLU is able to achieve 10 digits in the true relative accuracy, while
SlabLU achieves 8 digits. }
\label{fig:comparison_acc}
\end{figure}

SlabLU, on the other hand, uses an ordering based on a decomposition of the domain into slabs that has a limited pivoting scheme. Despite this limitation, SlabLU can achieve 10 digits of accuracy in the residual, which also gives high-order true relative accuracy, depending on the choice of $p$. SlabLU is a simple two-level framework that achieves large speedups over SuperLU by leveraging batched BLAS and GPU optimizations. Figure \ref{fig:comparison_build} provides a comparison between SuperLU and SlabLU in factorizing ${\tilde {\mtx A}}$ to solve (\ref{eq:hps_reduced}). Figure \ref{fig:comparison_acc} presents a comparison of accuracies in the computed solutions for various $p$.

\begin{figure}[htb!]
\centering
\includegraphics[width=0.46\textwidth]{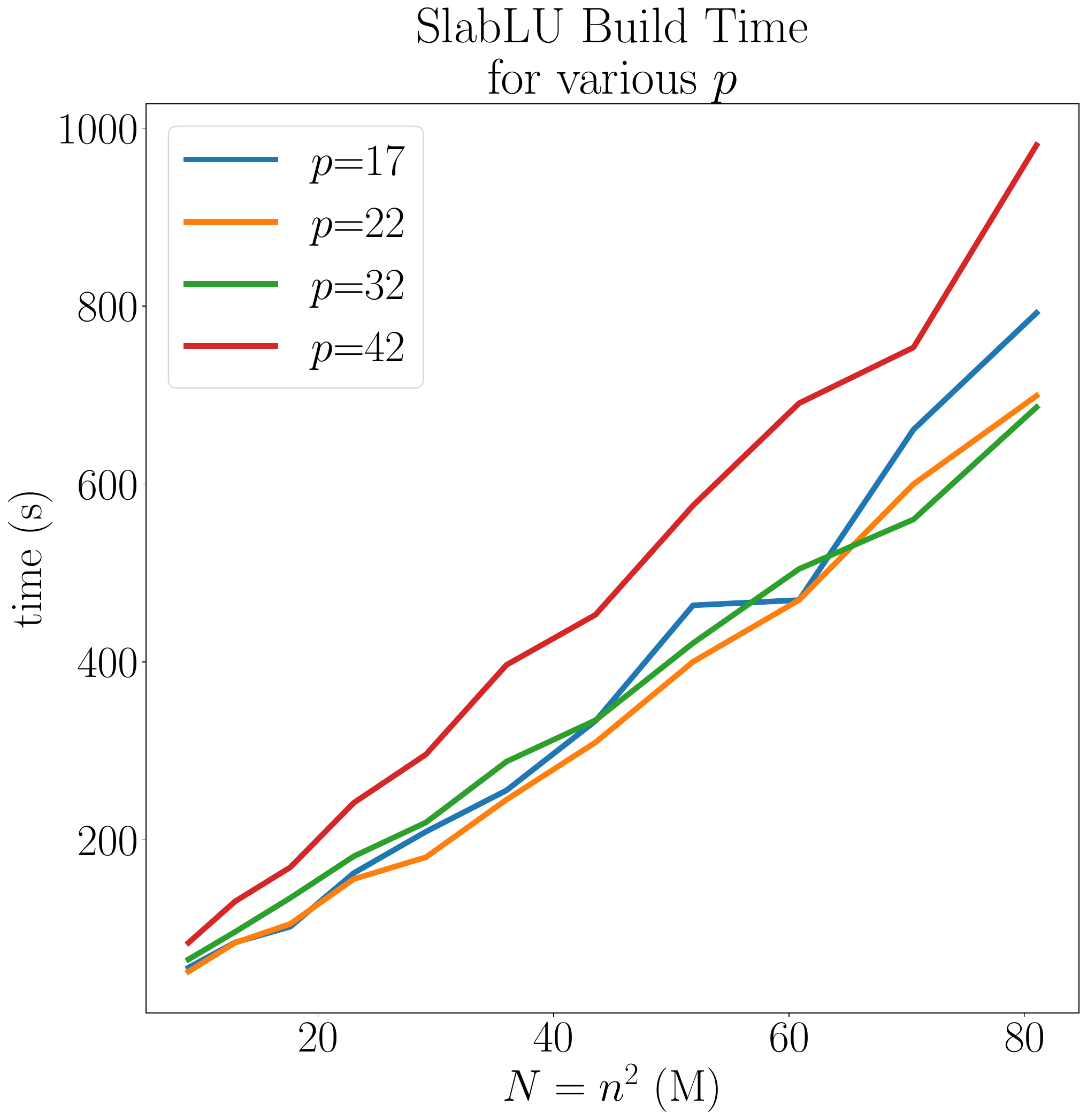}
\hfill
\includegraphics[width=0.49\textwidth]{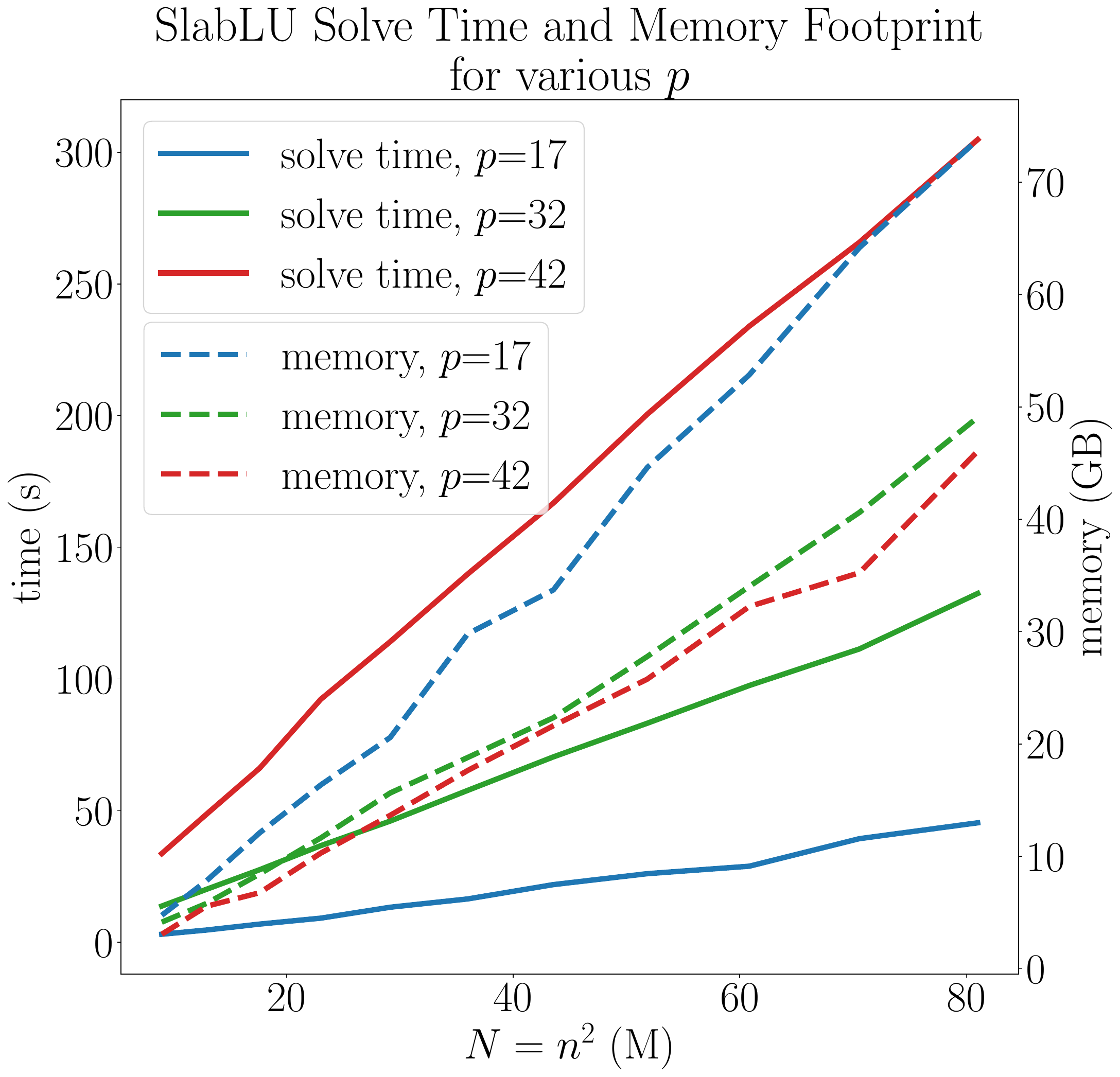}
\caption{Build time and solve time for HPS with various $p$ for (\ref{eq:constant_coeff_helm})
where $\kappa$ is increased with $N$ to maintain 10 points per wavelength.
The choice of $p$ does not substantially affect the time needed to factorize the 
sparse linear system with SlabLU. As $p$ increases, the memory footprint required to
store the factorization decreases and the solve time increases. 
}
\label{fig:slabLU_build}
\end{figure}

\begin{figure}[htb!]
\centering
\begin{minipage}{7.5cm}
\includegraphics[width=\textwidth]{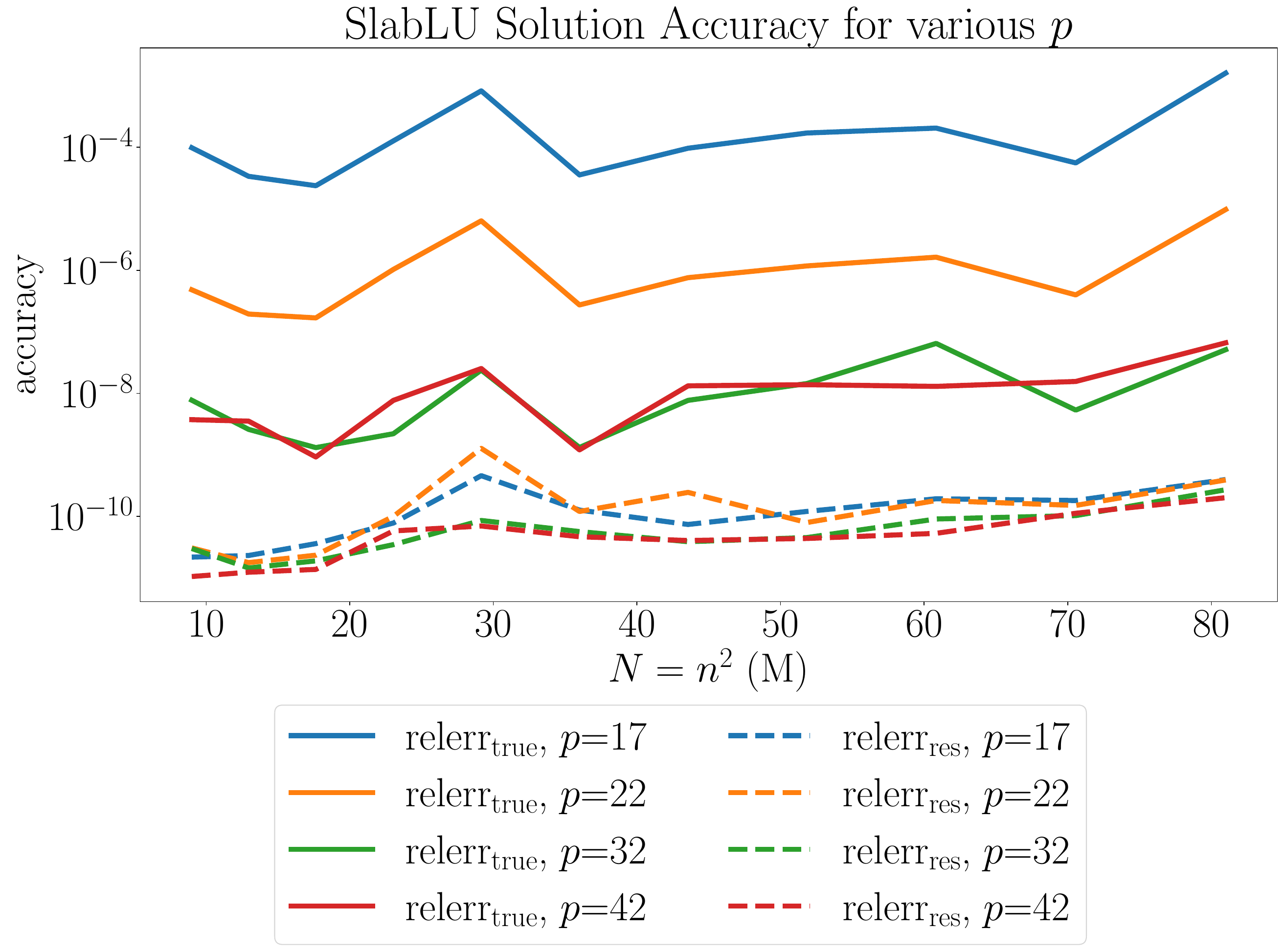}
\end{minipage}%
\hfill
\begin{minipage}{3.8cm}
\caption{Solution accuracy for HPS with various $p$ for (\ref{eq:constant_coeff_helm})
where $\kappa$ is increased with $N$ to maintain 10 points per wavelength.
Regardless of the choice of $p$, SlabLU resolves the solution to at least 10 digits of
relative accuracy in the residual. With increasing $p$, one can calculate solutions
with higher relative accuracy, compared to the true solution of the PDE.}
\label{fig:slabLU_acc}
\end{minipage}
\end{figure}

SlabLU can solve larger sparse systems (\ref{eq:hps_reduced}) with a smaller memory footprint than SuperLU. The memory footprint refers to how much main memory is required to store the sparse factorization of $\tilde {\mtx A}$. We demonstrate the ability of HPS, combined with SlabLU, for various $p$ to solve Helmholtz problems of size up to $900 \lambda \times 900 \lambda$ (for which $N$=81M) to high-order accuracy. Figure \ref{fig:slabLU_build} reports build and solve times for various choices of $p$, and Figure \ref{fig:slabLU_acc} reports the accuracy of the calculated solutions.


\subsection{Convergence for Scattering Problems for various $p$}

We will now demonstrate the ability of HPS, combined with SlabLU as a sparse direct solver, to solve 
complex scattering phenomena on various 2D domains. For the presented PDEs, we will show how the 
accuracy of the calculated solution converges to a reference solution depending on the choice of $p$ 
in the discretization. Specifically, we will solve the BVP (\ref{eq:bvp}) with the variable-coefficient 
Helmholtz operator (\ref{eq:var_helm}) for various Dirichlet data on smooth and rectangular domains.

\begin{figure}[htb!]
    \includegraphics[width=\textwidth]{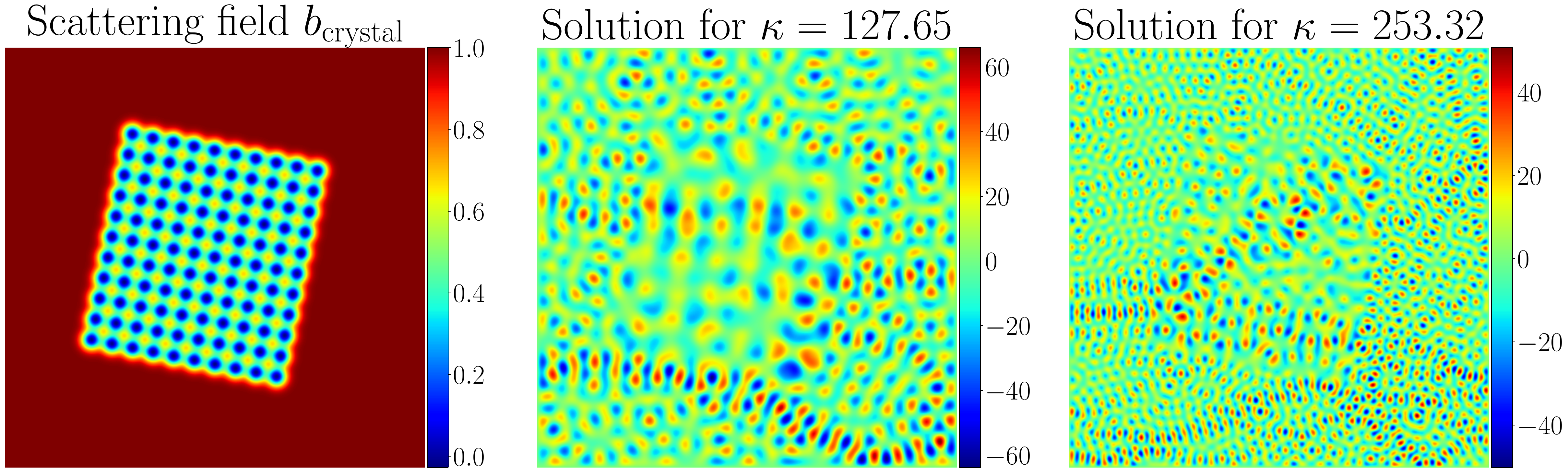}
    
    \vspace{0.5em}
    \begin{minipage}{0.5\textwidth}
    \centering
    \includegraphics[width=\textwidth]{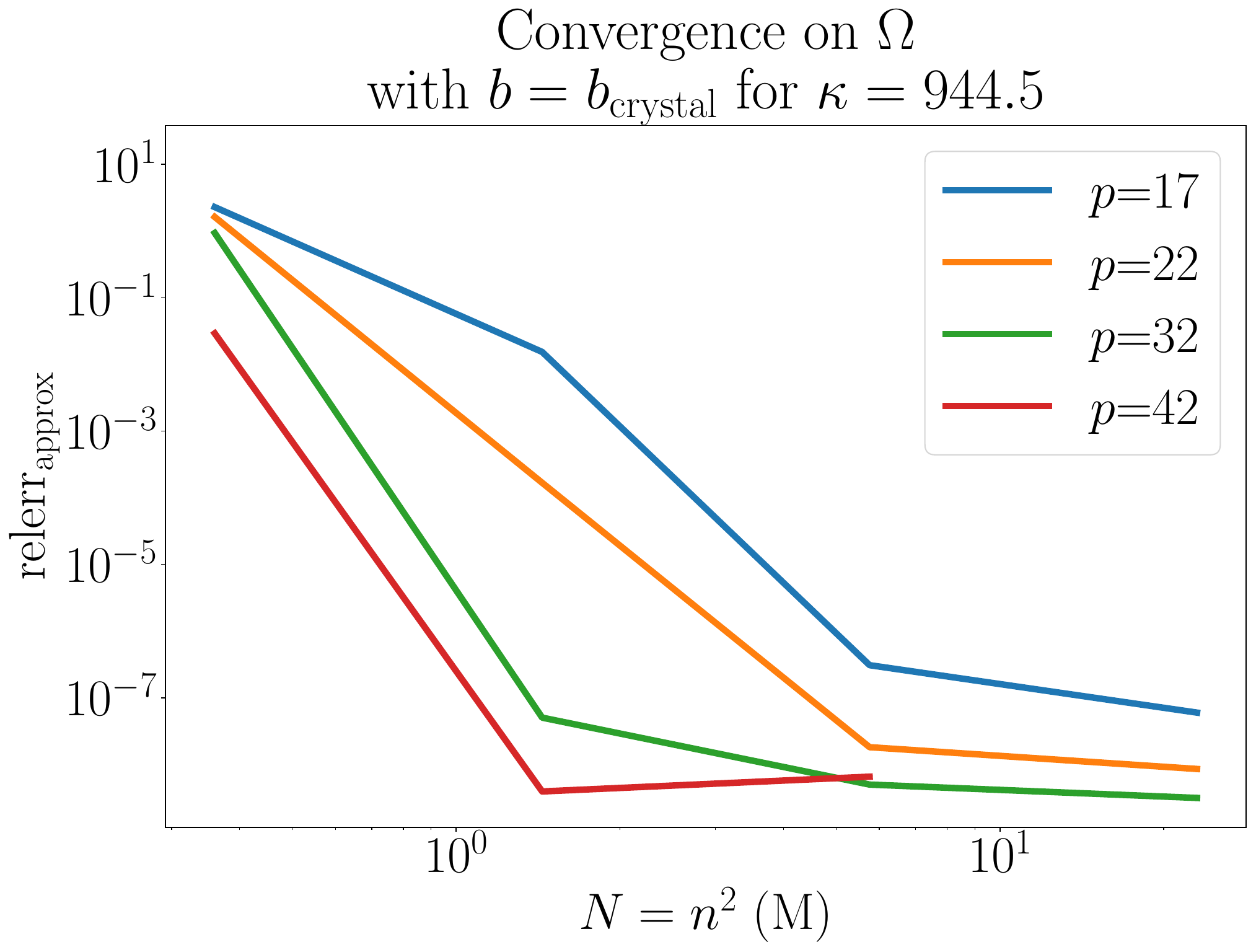}
    \end{minipage}%
    \hfill
    \begin{minipage}{0.48\textwidth}
    \caption{\small Solutions of variable-coefficient Helmholtz problem on square domain $\Omega$ with Dirichlet data given by $u \equiv 1$ on $\partial \Omega$ for various wavenumbers $\kappa$. The scattering field is $b_{\rm crystal}$, which is a photonic crystal that blocks the wave from propagating.}
    \label{fig:crystal_rhombus}
    \caption{Convergence on square domain $\Omega$ for reference solution $\mtx u_{\rm ref}$ on HPS discretization for $N$=36M with $p=42$.}
    \label{fig:convergence_square}
    \end{minipage}
\end{figure}

We fix the PDE and refining the mesh to compare calculated solutions to a reference solution obtained on a fine mesh with high $p$, as the exact solution is unknown. The relative error is calculated by comparing $\mtx u_{\rm calc}$ to the reference solution $\mtx u_{\rm ref}$ at a small number of collocation points $\{ x_j \}_{j=1}^M$ using the $l_2$ norm
\begin{equation}\rm {relerr}_{\rm approx} = \frac{{\| \mtx u_{\rm calc} - \mtx u_{\rm ref}\|}_2}{{\|\mtx u_{\rm ref}\|}_2}.\end{equation}
We demonstrate the convergence on a unit square domain $\Omega = [0,1]^2$ with a variable coefficient field $b_{\rm crystal}$ corresponding to a photonic crystal, shown in Figure \ref{fig:crystal_rhombus}. The convergence plot is presented in Figure \ref{fig:convergence_square}.

\begin{figure}[htb!]
\begin{minipage}{7.6cm}
\centering
\includegraphics[width=.965\textwidth]{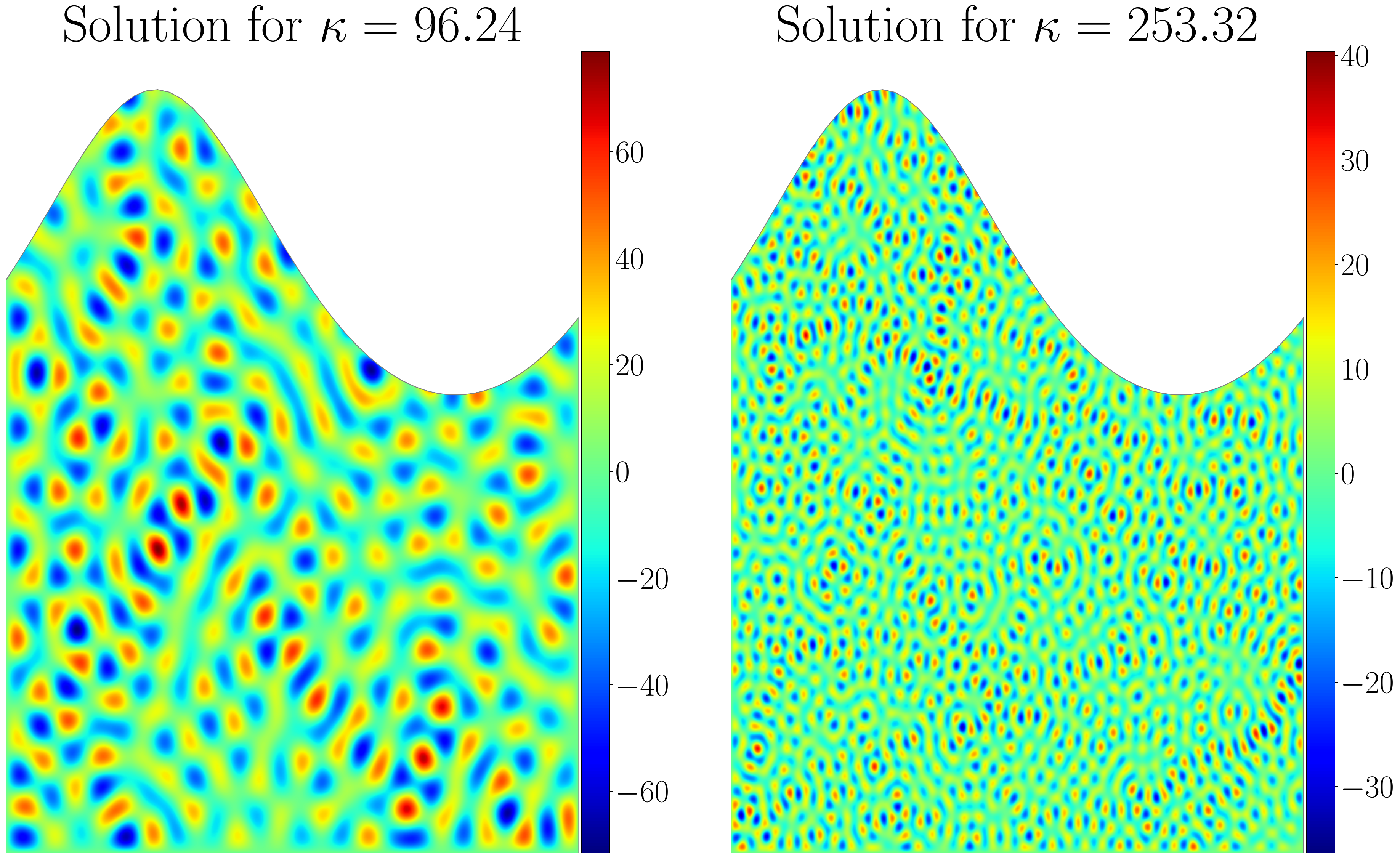}
\end{minipage}%
\hfill
\begin{minipage}{4.0cm}
\caption{\small Solutions of constant-coefficient Helmholtz problem on curved domain $\Psi$ with Dirichlet data given by $u \equiv 1$ on $\partial \Psi$ for various wavenumbers $\kappa$. The solutions are calculated parametrizing $\Psi$ in terms of a reference square domain $\Omega$ as (\ref{eq:parametrization}) and solving (\ref{eq:chain_rule}) on $\Omega$.}
\label{fig:sol_curved}
\end{minipage}

\vspace{1em}
\begin{minipage}{0.5\textwidth}
\centering
\includegraphics[width=\textwidth]{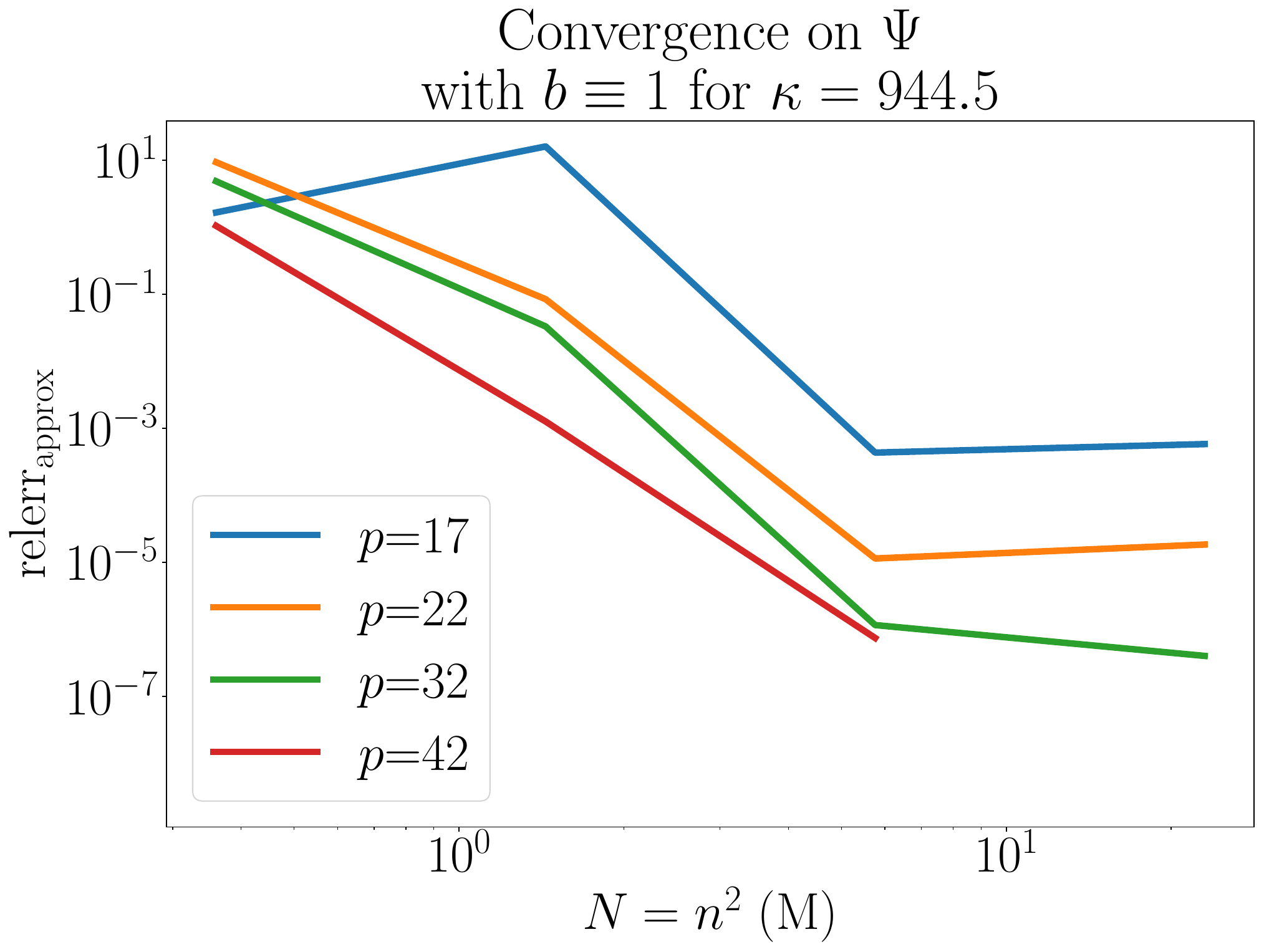}
\end{minipage}%
\hfill
\begin{minipage}{0.48\textwidth}
\caption{Convergence on curved domain $\Psi$ for reference solution $\mtx u_{\rm ref}$ on HPS discretization for $N$=36M with $p=42$. The solutions exhibit mildly singular behavior near the corners, and choosing high orders of $p$ aids in the convergence.}
\label{fig:convergence_curved}
\end{minipage}
\end{figure}

Next, we show the convergence on a curved domain $\Psi$ with a constant-coefficient field 
$b \equiv 1$, where $\Psi$ is given by an analytic parameterization over a reference square 
$\Omega = [0,1]^2$. The domain $\Psi$ is parametrized as
\begin{equation}
    \Psi =\left\{  \left(x_1, \frac{x_2}{\psi(x_1)} \right)\ \text{for}\ (x_1,x_2) \in \Omega = [0,1]^2 \right\},\ \text{where}\ \psi(z) = 1 - \frac 1 4 \sin(z).
\label{eq:parametrization}
\end{equation}
Using the chain rule, (\ref{eq:var_helm}) on $\Psi$ takes the following form on $\Omega$
\begin{align}
\begin{split}
    -\frac{\partial^2 u}{\partial x_1^2} &- 2 \frac{\psi'(x_1)x_2}{\psi(x_1)} \frac{\partial^2 u}{\partial x_1 \partial x_2}
    -  \left( \left( \frac{\psi'(x_1) x_2}{\psi(x_1)} \right)^2 + \psi(x_1)^2 \right) \frac{\partial^2 u}{\partial x_2^2} \\
    &- \frac{\psi''(x_1) x_2}{\psi(x_1)} \frac{\partial u}{\partial x_2} - \kappa^2 u = 0, \qquad (x_1,x_2) \in \Omega.
\end{split}
\label{eq:chain_rule}
\end{align}
The solutions on $\Psi$ are shown in Figure \ref{fig:sol_curved}, and the convergence plot is presented in Figure \ref{fig:convergence_curved}.

\section{Conclusions}

\vspace{-1em}

HPS is a high-order convergent discretization scheme that interfaces well with sparse direct solvers. In this manuscript, we describe GPU optimizations of the scheme that enable rapid and memory-efficient direct solutions of (\ref{eq:Au=f}) for resulting linear systems. First, we perform the leaf operations in parallel using batched BLAS. Then, we factorize a smaller system (\ref{eq:hps_reduced}) of size $\approx N/p$ using sparse direct solvers, where $p$ denotes the local order of convergence, which we show can be chosen as high as $42$.
The numerical results feature comparisons between sparse direct solvers and demonstrate that SlabLU can factorize systems corresponding to domains of size up to $900 \lambda \times 900 \lambda$ (for which $N$=81M) in less than 20 minutes. The approach is effective in resolving challenging scattering problems on various domains to high accuracy.

The techniques described are currently being implemented for three dimensional problems. The parallelizations described are immediately applicable. The scaling with $p$ deteriorates from $O(p^4 N)$ to $O(p^6 N)$, which limits how large $p$ can be chosen. However, initial numerical experiments demonstrate that $p=15$ remains viable on current hardware, which is high enough for most applications.

\vspace{0.5em}

\noindent
\textbf{Acknowledgments}
The work reported was supported by the Office of Naval Research (N00014-18-1-2354),
by the National Science Foundation (DMS-1952735 and DMS-2012606),
and by the Department of Energy ASCR (DE-SC0022251).
This is the author's accepted manuscript. The final published version is available at: \url{https://doi.org/10.1007/978-3-031-50769-4_62}.

\vspace{-2em}

\bibliographystyle{spphys}
\bibliography{main_bib}

\begin{thebibliography}{10}
\providecommand{\url}[1]{{#1}}
\providecommand{\urlprefix}{URL }
\expandafter\ifx\csname urlstyle\endcsname\relax
  \providecommand{\doi}[1]{DOI \discretionary{}{}{}#1}\else
  \providecommand{\doi}{DOI \discretionary{}{}{}\begingroup
  \urlstyle{rm}\Url}\fi

\bibitem{beriot2016efficient}
H.~B{\'e}riot, A.~Prinn, G.~Gabard, International Journal for Numerical Methods
  in Engineering \textbf{106}(3), 213 (2016)

\bibitem{deraemaeker1999dispersion}
A.~Deraemaeker, I.~Babu{\v{s}}ka, P.~Bouillard, International journal for
  numerical methods in engineering \textbf{46}(4), 471 (1999)

\bibitem{gillman2014direct}
A.~Gillman, P.G. Martinsson, SIAM Journal on Scientific Computing
  \textbf{36}(4), A2023 (2014)

\bibitem{martinsson2013direct}
P.G. Martinsson, Journal of Computational Physics \textbf{242}, 460 (2013)

\bibitem{ernst2012difficult}
O.G. Ernst, M.J. Gander, Numerical analysis of multiscale problems pp. 325--363
  (2012)

\bibitem{babb2018accelerated}
T.~Babb, A.~Gillman, S.~Hao, P.G. Martinsson, BIT Numerical Mathematics
  \textbf{58}, 851 (2018)

\bibitem{2019_martinsson_book}
P.G. Martinsson, \emph{Fast Direct Solvers for Elliptic PDEs}, \emph{CBMS-NSF
  conference series}, vol. CB96 (SIAM, 2019)

\bibitem{hao2016direct}
S.~Hao, P.G. Martinsson, Journal of Computational and Applied Mathematics
  \textbf{308}, 419 (2016)

\bibitem{cockburn2016static}
B.~Cockburn, in \emph{Building bridges: connections and challenges in modern
  approaches to numerical partial differential equations} (Springer, 2016), pp.
  129--177

\bibitem{guyan1965reduction}
R.J. Guyan, AIAA journal \textbf{3}(2), 380 (1965)

\bibitem{brubeck2021scalable}
P.D. Brubeck, P.E. Farrell, arXiv preprint arXiv:2107.14758  (2021)

\bibitem{2021_fortunato_ultraspherical}
D.~Fortunato, N.~Hale, A.~Townsend, Journal of Computational Physics
  \textbf{436}, 110087 (2021)

\bibitem{dongarra2017design}
J.~Dongarra, S.~Hammarling, N.J. Higham, S.D. Relton, P.~Valero-Lara,
  M.~Zounon, Procedia Computer Science \textbf{108}, 495 (2017)

\bibitem{2006_davis_directsolverbook}
T.A. Davis, \emph{Direct methods for sparse linear systems}, vol.~2 (Siam,
  2006)

\bibitem{2016_acta_sparse_direct_survey}
T.A. Davis, S.~Rajamanickam, W.M. Sid-Lakhdar, Acta Numerica \textbf{25}, 383
  (2016).
\newblock \doi{10.1017/S0962492916000076}

\bibitem{yesypenko2022slablu}
A.~Yesypenko, P.G. Martinsson, arXiv preprint arXiv:2211.07572  (2022)

\bibitem{levitt2022linear}
J.~Levitt, P.G. Martinsson, arXiv preprint arXiv:2205.02990  (2022)

\bibitem{2011_li_supernodal}
X.S. Li, M.~Shao, ACM Transactions on Mathematical Software (TOMS)
  \textbf{37}(4), 1 (2011)

\end{thebibliography}
\end{document}